\documentclass[12pt,leqno]{article}
\usepackage{amsmath, amsthm, amsfonts, amssymb, color}
\usepackage{mathrsfs}
\theoremstyle{definition}

\newcommand{\scr}[1]{\mathscr #1}
\definecolor{wco}{rgb}{0.5,0.2,0.3}

\numberwithin{equation}{section} \theoremstyle{remark}

\newcommand{\ua}{\uparrow}

\title{{\bf Transportation-Cost Inequalities on  Path Space Over Manifolds with Boundary}\footnote{Supported in
 part by WIMICs, NNSFC(10721091) and the 973-Project.}
}
\author{
{\bf Feng-Yu Wang}\\
\footnotesize{School of Mathematical Sci. and Lab. Math. Com. Sys.,
Beijing Normal
University, Beijing 100875, China}\\
\footnotesize{and}\\ \footnotesize{Department of Mathematics,
Swansea University, Singleton Park, SA2 8PP, UK}\\ \footnotesize{Email: wangfy@bnu.edu.cn;
F.Y.Wang@swansea.ac.uk}}
\begin{document}
\def\R{\mathbb R}  \def\ff{\frac} \def\ss{\sqrt} \def\B{\mathbf
B}
\def\N{\mathbb N} \def\kk{\kappa} \def\m{{\bf m}}
\def\dd{\delta} \def\DD{\Delta} \def\vv{\varepsilon} \def\rr{\rho}
\def\<{\langle} \def\>{\rangle} \def\GG{\Gamma} \def\gg{\gamma}
  \def\nn{\nabla} \def\pp{\partial} \def\tt{\tilde}
\def\d{\text{\rm{d}}} \def\bb{\beta} \def\aa{\alpha} \def\D{\scr D}
\def\EE{\mathbb E} \def\si{\sigma} \def\ess{\text{\rm{ess}}}
\def\beg{\begin} \def\beq{\begin{equation}}  \def\F{\scr F}
\def\Ric{\text{\rm{Ric}}} \def\Hess{\text{\rm{Hess}}}
\def\e{\text{\rm{e}}} \def\ua{\underline a} \def\OO{\Omega}  \def\oo{\omega}
 \def\tt{\tilde} \def\Ric{\text{\rm{Ric}}}
\def\cut{\text{\rm{cut}}} \def\P{\mathbb P} \def\ifn{I_n(f^{\bigotimes n})}
\def\C{\scr C}      \def\aaa{\mathbf{r}}     \def\r{r}
\def\gap{\text{\rm{gap}}} \def\prr{\pi_{{\bf m},\varrho}}  \def\r{\mathbf r}
\def\Z{\mathbb Z} \def\vrr{\varrho} \def\ll{\lambda}
\def\L{\scr L}\def\Tt{\tt} \def\TT{\tt}\def\II{\mathbb I}
\def\i{{\rm i}}\def\Sect{{\rm Sect}}\def\E{\mathbb E} \def\H{\mathbb H}
\def\M{\scr M}\def\Q{\mathbb Q} \def\texto{\text{o}}

\maketitle
\begin{abstract} Let $L=\DD+Z$ for a $C^1$ vector field $Z$ on a complete
Riemannian manifold  possibly with a boundary. By using the uniform
distance, a number of  transportation-cost  inequalities on the path
space for the (reflecting) $L$-diffusion process are proved to be
equivalent to the curvature condition $\Ric-\nn Z\ge - K$ and the
convexity of the boundary (if exists). These inequalities are new
even for manifolds without boundary, and are partly extended to
non-convex manifolds by using a conformal change of metric which
makes the boundary from non-convex to convex.
\end{abstract} \noindent
 \vskip 2cm

\section{Introduction}

In 1996  Talagrand   \cite{T} found that the $L^2$-Wasserstein distance to the standard Guassian measure can be
 dominated by the square root of twice relative entropy. This
inequality is called (Talagrand) transportation-cost inequality, and
has been extended to distributions on finite- and
infinite-dimensional spaces. In particular, this inequality was
established on the path space of diffusion processes with respect to
several different distances (i.e. cost functions): see e.g. \cite{F}
for the study on the Wiener space with the Cameron-Martin distance,
\cite{W02, DGW} on the path space of diffusions with the
$L^2$-distance, \cite{W04} on the Riemannian
 path space with intrinsic distance induced by the Malliavin gradient operator, and \cite{FWW,
WZ} on the path space of diffusions with the uniform distance. The
main purpose of this paper is to investigate the Talagrand
inequality on the path space of reflecting diffusion process, for
which both the curvature and the second fundamental form of the
boundary will take important roles.

Let $M$ be a connected complete Riemannian manifold possibly with a
boundary  $\pp M$. Let $L=\DD +Z$ for a $C^1$ vector field $Z$ on
$M$. Let $X_t$ be the (reflecting if $\pp M\ne\emptyset$) diffusion
process generated by $L$ with initial distribution $\mu\in \scr
P(M),$ where $\scr P(M)$ is the set of all probability measures on
$M$. Assume that $X_t$ is non-explosive, which is the case if $\pp
M$ is convex and the curvature condition
 \beq\label{C} \Ric-\nn Z\ge -K\end{equation} holds for some constant $K\in
 \R.$ In this case, for any $T>0$, the distribution  $\Pi_\mu^T$ of $X_{[0,T]}:= \{X_t:\ t\in [0,T]\}$  is a
 probability measure on the (free) path space

$$M^T:=C([0,T];M).$$
 When $\mu=\dd_o$, the Dirac measure at point $o\in M$, we
 simply denote $\Pi_{\dd_o}^T=\Pi_o^T.$ For any
 nonnegative measurable function $F$ on $M_T$ such that
 $\Pi_\mu^T(F)=1,$ one has

 \beq\label{0} \mu_F^T(\d x):= \Pi_x^T(F)\mu(\d x)\in \scr
 P(M).\end{equation}

 Let $\rr$ be the Riemannian distance on $M$; i.e. for $x,y\in M, \rr(x,y)$ is the length
 of the shortest curve on $M$ linking $x$ and $y$. Then $M^T$ is a Polish space
 under the uniform distance

 $$\rr_\infty(\gg,\eta)= \sup_{t\in [0,T]} \rr(\gg_t,\eta_t),\ \ \ \gg,\eta\in M^T.$$
 Let $W_{2,\rr_\infty}$ be the $L^2$-Wasserstein distance (or $L^2$-transportation cost) induced by
 $\rr_\infty$. In general, for any $p\ge 1$ and
 for two   probability measures $\Pi_1,\Pi_2$ on $M^T$,

 $$W_{p,\rr_\infty}(\Pi_1,\Pi_2):=
 \inf_{\pi\in \scr C(\Pi_1,\Pi_2)} \bigg\{\int_{M_{T}\times M_{T}}
 \rr_\infty(\gg,\eta)^p\pi(\d\gg,\d\eta)\bigg\}^{1/p}$$
 is the $L^p$-Warsserstein distance (or $L^p$-transportation cost) of $\Pi_1$ and $\Pi_2$ induced
 by the uniform norm,
 where $\scr C(\Pi_1,\Pi_2)$ is the set of all couplings for $\Pi_1$ and $\Pi_2$.

Before moving on, let us recall the Talagrand transportation-cost
inequality established in \cite{FWW} on the path space over
Riemannian manifolds without boundary.  Let $\pp M=\emptyset$ and
$\rr_o=\rr(o,\cdot).$ If

 \beq\label{A} |Z|\le \psi\circ\rr_o\end{equation} holds for some positive function $\psi$ such that
  $\int_0^\infty \ff 1 {\psi(s)}\d s=\infty$, then (see \cite[Theorem
  1.1]{FWW})
  \beq\label{1.3}  W_{2, \rr_\infty}(F\Pi_o^T,\Pi_o^T)^2\le \ff 2 K (\e^{2KT}-1)\Pi_o^T(F\log F),
  \ \ F\ge 0, \Pi_o^T(F)=1.\end{equation}

According to \cite{OV, BGL, W04},     the log-Sobolev inequality for
a smooth elliptic diffusion implies the   Talagrand
transportation-cost inequality with the intrinsic distance. So,
(\ref{1.3}) was proved in \cite{FWW} by using a known damped
log-Sobolev inequality on the path space and finite-dimensional
approximations. To ensure the smoothness of the approximating
diffusions, one needs
  the boundedness of curvature. To get rid of this
condition, a sequence of new metric approximating the original one
were constructed in \cite{FWW},  which  satisfy  (\ref{C}) and have
bounded curvatures. In this way  (\ref{1.3}) was established without
using curvature upper bounds. But to realize this approximation
argument,   the technical condition (\ref{A}) with $\int_0^\infty
\ff 1 {\psi(s)}\d s=\infty$ was adopted.

In this paper we adopt a different argument developed in \cite{WZ}
for diffusions on $\R^d$ by using the martingale representation
theorem and Girsanov transformations, so that this technical
condition was avoided. Furthermore, we present a number of cost
inequalities which are equivalent to the convexity of $\pp M$ (if
exists) and the curvature condition (\ref{C}).

When $\pp M\ne\emptyset$, let $N$ be the inward unit normal vector
field of $\pp M$. Then the second fundamental form of $\pp M$ is
defined by

$$\II(U,V)=-\<\nn_U N,V\>,\ \ \ U,V\in T\pp M,$$ where $T\pp M$ is the tangent space of $\pp M.$
If $\II\ge 0$, i.e. $\II(U,U)\ge 0$ for all $U\in T\pp M$,  we call
$M$ (or $\pp M$) convex.

\beg{thm} \label{T1.1} Let $P_T(o,\cdot)$ be the distribution of
$X_T$ with $X_0=o$, and let $P_T$ be the corresponding semigroup.
The following statements are equivalent to each other:
\beg{enumerate}
\item[$(1)$]\ $\pp M$ is either convex or empty, and  $(\ref{C})$ holds.
\item[$(2)$]\ For any $T>0, \mu\in\scr P(M)$ and nonnegative $F$
with $\Pi_\mu^T(F)=1$,
$$W_{2, \rr_\infty}(F\Pi_\mu^T,\Pi_{\mu_F^T}^T)^2\le \ff 2 K (\e^{2KT}-1)\Pi_\mu^T(F\log
F)$$ holds, where $\mu_F^T\in\scr P(M)$ is fixed by $(\ref{0})$.
\item[$(3)$]\   $(\ref{1.3})$ holds for any $o\in M$ and $T>0.$
\item[$(4)$]\ For any $o\in M$ and $T>0$,

$$W_{2,\rr}\big(P_T(o,\cdot), fP_T(o,\cdot)\big)^2\le \ff 2 K (\e^{2KT}-1)P_T(f\log f)(o),
\ \ \ f\ge 0, P_Tf(o)=1.$$
\item[$(5)$]\   For any $T>0$,  $\mu,\nu\in\scr P(M),$ and $p\ge 1$,

$$W_{p, \rr_\infty} (\Pi_\mu^T,\Pi_\nu^T) \le \e^{KT} W_{p,\rr}(\mu,\nu),$$
where $W_{p,\rr}$ is the $L^p$-Wasserstein distance for probability
measures on $M$ induced by $\rr$.
\item[$(6)$]\ For any $x,y\in M$ and   $T>0$,

$$W_{2,\rr}\big(P_T(x,\cdot), P_T(y,\cdot)\big)\le \e^{KT}\rr(x,y).$$
\item[$(7)$] For any $T>0$, $\mu\in\scr P(M),$ and $F\ge 0$ with
$\Pi_\mu^T(F)=1,$

$$W_{2,\rr_\infty} (F\Pi_\mu^T, \Pi_\mu^T)\le\Big\{\ff 2 K (\e^{2KT}-1)\Pi_\mu^T(F\log
F)\Big\}^{1/2} + \e^{KT} W_{2,\rr}(\mu_F^T,\mu).$$
\item[$(8)$] For any $\mu\in \scr P(M)$ and $C\ge 0$ such that

$$ W_{2,\rr}(f\mu, \mu)^2\le C\mu(f\log f),\ \ f\ge 0,
\mu(f)=1,$$ there holds

$$W_{2,\rr_\infty}(F\Pi_\mu^T, \Pi_\mu^T)^2\le \bigg(\ss{\ff 2 K
(\e^{2KT}-1)} +\ss{C}\bigg)^2 \Pi_\mu^T(F\log F),\ \ F\ge 0,
\Pi_\mu^T(F)=1.$$
\end{enumerate} \end{thm}

When $\pp M=\emptyset,$ there exist many equivalent semigroup
inequalities for the curvature condition (\ref{C}):   see e.g.
\cite{B,L} for equivalent statements on gradient estimates,
log-Sobolev/Poicar\'e inequalities, and isoperimetric inequality;
   \cite{W04b, W09} for equivalent Harnack type inequalities;
and   \cite{RS} for equivalent inequalities on Wasserstein
distances. Theorem \ref{T1.1}  provides seven   equivalent
inequalities for the convexity of $\pp M$ (if exists) and the
curvature condition (\ref{C}), which are new even for manifolds
without boundary.

To prove this   Theorem, we shall use a formula of the second
fundamental form established in \cite{W09} for compact manifolds
with boundary. Since in this paper the manifold is allowed to be
non-compact, we shall reprove this formula in Section 2 by using the
reflecting diffusion process up to the exit time of a compact
domain. This formula implies the equivalence of Theorem
\ref{T1.1}(1) and the semigroup log-Sobolev/Poincar\'e inequalities
(see Theorem \ref{T2.4} below). In Section 3 we prove Theorem
\ref{T1.1} by using results in Section 2, the martingale
representation and Girsanov transformation for (reflecting)
  diffusions on (convex) manifolds.
 which lead to a proof from (1) to (2), then prove (1) from (4) by using results obtained in
 Section 2. The proof of Theorem \ref{T1.2} will be
addressed in Section 4.

To establish transportation-cost inequalities on the path space for
non-convex manifolds, we shall adopt a conformal change of metric
$\<\cdot,\cdot\>'= f^{-2}\<\cdot,\cdot\>$ such that $\pp M$ is
convex under the new metric (see \cite[Lemma 2.1]{W07}). Let $\DD'$
be the Laplacian induced by the new metric, we have (see \cite[Lemma
2.2]{W07})

\beq\label{L} L= f^{-2} \Big\{\DD' + \varphi^2Z + \ff{d-2}2 \nn
f^2\Big\}.\end{equation} Thus, in Section 4  we modify our arguments
 to study the reflecting diffusion process with a non-constant coefficient, from
 which we  partly extend  Theorem \ref{T1.1}  to non-convex manifolds in Section
 5
 to non-convex manifolds.

 \section{Formulae for the second fundamental form and applications}

When $M$ is compact, the following formula on $\pp M$ has been found
 in \cite{W09}:

\beq\label{*2} \lim_{t\to 0} \ff{|\nn f|^2}{\ss t} \log \ff{|\nn
P_t|}{(P_t |\nn f|^p)^{1/p}}=-\ff 2 {\ss\pi} \II(\nn f,\nn f),\ \
p\ge 1,\end{equation} where $f$ is a smooth function satisfying the
Neumann boundary condition. When $M$ is non-compact, some  technical
problems appear in the original proof when e.g. a dominated
convergence is used. To fix these problems, we shall stop the
process in a compact domain, so that we shall first study the
behavior of hitting times.

Recall that the reflecting $L$-diffusion process can be constructed
by solving the SDE

\beq\label{2.0} \d X_t =\ss 2\, \Phi_t\circ \d B_t +  Z(X_t)\d t
+N(X_t)\d l_t,
\end{equation} where $\Phi_t$ is the horizontal lift of $X_t$ onto the frame bundle
$O(M)$,  $B_t$ is the $d$-dimensional Brownian motion.

By the It\^o formula, for any $f\in C^2(M)$ we have

\beq\label{I} \d f(X_t) = \ss 2\<\nn f(X_t),\Phi_t\circ \d B_t\> +
Lf(X_t)\d t + Nf(X_t)\d l_t,\end{equation} where $Nf = \<N,\nn f\>$.
For any $R>0$, let

$$\tau_R=\inf\{t\ge 0:\ \rr(X_0,X_t)\ge R\}.$$

\beg{prp}\label{P1} Let $R>0$ and $X_0=o\in M$ be fixed. Then there
exist two constants $c_1,c_2>0$ such that

$$\P(\tau_R\le t)\le c_1\e^{-c_2/t},\ \ \ t>0.$$\end{prp}

\beg{proof} This result is well known on manifolds without boundary
(cf. \cite[Lemma 2.3]{ATW09}), and the proof works also when $\pp M$
is convex. As in the present case the boundary is not necessarily
convex, we shall follow \cite{W07} to make the boundary convex under
a conformal change of metric. Since

$$\B_R:=\{x\in M:\ \rr(o,x)\le R\}$$ is compact, there exists a
constant $\si>0$ such that $\II\ge -\si$ holds on $\pp M\cap \B_R.$
Let $f\ge 1$ be smooth such that

\beq\label{A1} N\log f\ge \si\ \ \text{on}\ \pp M\cap
\B_R.\end{equation} Such a function can be constructed by using the
distance function $\rr_\pp$ to the boundary $\pp M$. Since $\B_{2R}$
is compact, there exists a constant  $r_0>0$ such that $\rr_\pp$ is
smooth on $\{x\in \B_{2R}:\ \rr_\pp(x)\le r_0\}.$ Let $h\in
C^\infty([0,\infty))$ such that $h'\ge 0, h(0)=1, h'(0)=\si$ and
$h'(r)=0$ for $r\ge r_0.$ Then $h\circ \rr_\pp$ is smooth on
$\B_{2R}$ and $N\log h\circ\rr_\pp|_{\pp M\cap \B_{2R}}= \si$. Thus,
it suffices to take smooth $f\ge 1$ such that $f= h\circ\rr_\pp$ on
$\B_R.$

By \cite[Lemma 2.1]{W07} and (\ref{A1}), $\pp M$ is convex in $\B_R$
under the new metric

$$\<\cdot,\cdot\>':= f^{-2}\<\cdot,\cdot\>,$$ where $\<\cdot,\cdot\>$
is the original metric. Let $\DD'$ be the Laplacian induced by the
new metric. We have (see \cite[Lemma 2.2]{W07})

$$L=f^{-2} (\DD'+Z')$$ for some $C^1$-vector field $Z'.$ Let
$\tt\rr_o$ be the Riemannian distance to $o$  induced by the new
metric. By the Laplacian comparison theorem,

\beq\label{A2} L\tt\rr_o^2\le c\ \ \text{on}\ \B_R\end{equation}
holds for some constant $c>0$ outside the cut-locus induced by
$\<\cdot,\cdot\>'.$ Since $\pp M$ is convex on $\B_R$ and $N$ is
still the inward normal vector under the new metric, we have

$$N\tt\rr_o\le 0\ \ \text{on}\ \pp M\cap \B_R.$$Therefore, by using
Kendall's It\^o formula for the distance (cf. \cite{K} for $f=1$),
(\ref{A2}) implies

$$\d\tt\rr_o^2(X_t)\le 2\ss 2\, f^{-2}(X_t) \tt\rr_o(X_t)\d b_t+c\d
t,\ \ \ t\le \tau_R,$$ where $b_t$ is some one-dimensional Brownian
motion. Since  $f^{-2}\le 1$, this implies that for any $\dd>0$, the
process

$$Z_s:= \exp\bigg[\ff \dd t \tt\rr_o^2(X_s) -\ff\dd t cs -4
\ff{\dd^2}{t^2} \int_0^s\tt\rr_o^2(X_u)\d u\bigg],\ \ s\le \tau_R$$
is a super martingale. Therefore, letting $C>1$ be a constant such
that $f\le C$ on $\B_R$ and thus, $\rr_o\ge\tt\rr_o\ge C^{-1}\rr_o$
holds on $\B_R$, we obtain

\beg{equation*}\beg{split} \P(\tau_R\le t) &=
 \P\Big(\max_{s\in [0,t]}\rr_o(X_{s\land\tau_R})\ge R\Big)\le
 \P\Big(R\ge \max_{s\in [0,t]}\tt\rr_o(X_{s\land\tau_R})\ge \ff R
 C\Big)\\
 &\le \P\Big(\max_{s\in [0,t]} Z_{s\land \tau_R} \ge
 \exp\Big[\ff{\dd R^2}{tC^2} -\dd c
 -\ff{4\dd^2R^2}{t}\Big]\Big)\\
 &\le \exp\Big[ c\dd -\ff {R^2} {tC^2} (\dd - 4C^2\dd^2)\Big],\ \ \
 \dd>0.\end{split}\end{equation*} The proof is then completed by
 taking e.g. $\dd= 1/(8C^2).$\end{proof}

 \beg{prp}\label{P2} Let $X_0=o\in \pp M$. Then for any $R>0$,

$$\limsup_{t\to 0} \ff 1 t \big|\E l_{t\land \tau_R}
-2\ss{t/\pi}\big|<\infty.$$\end{prp}

\beg{proof} Repeating the proof of \cite[Lemma 2.2]{W09} by using
$t\land \tau_R$ in place of $t$, we obtain

\beq\label{A4} \E l_{t\land \tau_R}^2\le ct,\ \ \ t\in
[0,1]\end{equation}  for some constant $c>0.$ Let $r_0>0$ be such
that $\rr_\pp$ is smooth on $\{\rr_\pp\le r_0\}\cap \B_R$. Let

$$\tau=\inf\{t\ge 0:\ \rr_\pp(X_t)\ge r_0\}.$$ By the It\^o formula
we have

\beq\label{A5} \d\rr_\pp(X_t)=\ss 2 \,\d b_t +L\rr_\pp(X_t)\d t +\d
l_t,\ \ \ t\le \tau\land \tau_R,\end{equation} where, as before,
$b_t$ is some one-dimensional Brownian motion. By the proof of
\cite[Theorem 2.1]{W09} using $\tau\land\tau_R$ in place of $\tau$,
we have, instead of (2.4) in \cite{W09},

\beq\label{A6} \E\big(\rr_\pp (X_{t\land \tau\land \tau_R})-\ss 2\,
|\tt b_{t\land \tau\land \tau_R}|\big)^2 \le c_1t^2,\ \ t\in
[0,1]\end{equation} for some constant $c_1>0,$ where $\tt b_t$ is
some one-dimensional Brownian motion. Due to (\ref{A5}),

$$\big|\E l_{t\land\tau\land\tau_R} -\E \rr_\pp
(X_{t\land\tau\land\tau_R})\big|\le c_2 t$$ holds for some constant
$c_2>0$. Combining this with (\ref{A6}) we arrive at

$$\big|\E l_{t\land\tau\land\tau_R} -\ss 2\, \E |\tt b_{t\land \tau\land \tau_R}|\big|\le c_3
t,\ \ \ t\in [0,1]$$ for some constant $c_3>0.$ Since $\E|\tt b_t|=
\ss{2t/\pi}$ and $\E |\tt b_t|^2=t$, this and (\ref{A4}) imply

\beq\label{A7} \beg{split} &\Big|\E l_{t\land\tau_R} - \ff{2\ss
t}{\ss\pi}\Big| = \big|\E l_{t\land \tau_R} - \ss 2\, \E|\tt
b_t|\big|\\
&\le c_3 t + \E 1_{\{t\ge \tau\land\tau_R\}} (l_{t\land \tau_R}+\ss
2\,
|\tt b_{t}|)\\
&\le c_3t +c_4 \ss{t\P(t\ge \tau\land \tau_R)},\ \ \ t\in
[0,1].\end{split}\end{equation}  Moreover, noting that

$$\P(\tau\land \tau_R\le t, \tau_R>\tau)\le \P\Big(\max_{s\in [0,t]} \rr_\pp
(X_{s\land\tau\land\tau_R})\ge r_0\Big),$$ by using $\tau\land
\tau_R$ to replace $\tau$ in the proof of \cite[Proposition
A.2]{W09}, we conclude that

$$\P(\tau\land\tau_R\le t,\ \tau_R> \tau)\le c_5
\exp[-r_0^2/(16t)],\ \ \ t>0$$ holds for some constant $c_5>0.$
Combining this with Proposition \ref{P1}, we obtain

$$\P(t\ge \tau\land\tau_R)\le c_6\e^{-c_7/t},\ \ \ t>0$$ for some
constants $c_6,c_7>0.$ Therefore, the proof is completed by
(\ref{A7}).\end{proof}

\beg{thm}\label{T} Let $f\in C^\infty(M)$ with $Nf|_{\pp M}=0$.

$(1)$ For any $p\ge 1$ and $R>0$,

\beq\label{*1} \lim_{t\to 0}\ff{|\nn f|^2}{\ss t} \log \ff{(\E|\nn
f|^p(X_{t\land\tau_R})|)^{1/p}}{|\nn f|}=\ff 2 {\ss\pi} \II(\nn
f,\nn f)\end{equation} holds at points on $\pp M$ such that $|\nn
f|>0.$

$(2)$ Assume that for any $g\in C_0^1(M)$ the function $|\nn P_\cdot
g|$ is bounded on $[0,1]\times M$. If moreover $f$ has a compact
support, then $(\ref{*2})$ holds  points on $\pp M$ such that $|\nn
f|>0.$\end{thm}

\beg{proof} (\ref{*1}) follows immediately from the proof of
\cite[Theorem 1.2]{W09} by using Proposition \ref{P2} in place of
\cite[Theorem 2.1]{W09}, and using $t\land \tau_R$ in place of $t$.

Next, let $f\in C_0^\infty(M)$.  By the assumption of (2) and that
$Lf\in C_0^1(M)$, $|\nn P_\cdot Lf|$ is bounded on $[0,1]\times M$.
So, the proof of \cite[(3.1)]{W09} implies that

\beq\label{A8} \lim_{t\to 0} \ff{|\nn f|^2}{\ss t} \log \ff{|\nn P_t
f|}{(P_t |\nn f|^p)^{1/p}} =-\lim_{t\to 0} \ff{|\nn f|^2}{\ss t}\log
\ff{(P_t |\nn f|^p)^{1/p}}{|\nn f|}.\end{equation} Since by
Proposition \ref{P1}, there exist two constant $c_1,c_2>0$ such that

$$\big|P_t |\nn f|^p -\E|\nn f|^p(X_{t\land\tau_R})\big|\le \|\nn
f\|_\infty^p \P(t>\tau_R)\le c_1\e^{-c_2/t},\ \ t>0,$$ we conclude
that (\ref{*2}) follows from (\ref{A8}) and (\ref{*1}).\end{proof}

As an application of (\ref{*1}), the following result provides
equivalent semigroup log-Sobolev/Poincar\'e inequalities for Theorem
\ref{T1.1}(1).

\beg{thm}\label{T2.4} Each of the following statements is equivalent
to Theorem $\ref{T1.1}(1)$:\beg{enumerate} \item[$(9)$] For any
$T>0$ and $f\in C_b(M)$,

$$P_T f^2\log f^2\le (P_T f^2)\log P_T f^2 +\ff {\e^{2KT}-1}{2K}
P_T|\nn f|^2.$$
\item[$(10)$] For any
$T>0$ and $f\in C_b(M)$,

$$P_T f^2 \le (P_T f)^2 +\ff {\e^{2KT}-1}{K}
P_T|\nn f|^2.$$\end{enumerate}\end{thm}

\beg{proof} According to e.g. \cite[Lemma 3.1]{W97}, which holds
also for the non-symmetric case, Theorem \ref{T1.1}(1) implies the
semigroup log-Sobolev inequality (9). It is well known that the
log-Sobolev inequality implies the Poincar\'e inequality. So,  (10)
follows from (9). Hence, it remains to show that (10) implies
Theorem \ref{T1.1}(1). Below we shall prove the convexity of $\pp M$
and the curvature condition (\ref{C}) respectively.

 (a) Let $\pp M\ne\emptyset.$
For any $o\in \pp M$ and non-trivial $X\in T_o\pp M$, we aim to show
that $\II(X,X)\ge 0.$ Let $f\in C_b^\infty(M)$ such that $Nf|_{\pp
M}=0$ and $\nn f(o)=X.$ Let $X_0=o$ and

$$\tau_1=\inf\{t\ge 0: \rr(o,X_t)\ge 1\}.$$ Since $f$ and $f^2$ satisfies
the Neumann boundary condition, we have

\beg{equation*}\beg{split} &\E f(X_{t\land \tau_1})= f(o)
+\E\int_0^{t\land \tau_1} Lf(X_s)\d s,\\
&\E f^2(X_{t\land \tau_1})= f^2(o) +2\E\int_0^{t\land \tau_1}
(fLf)(X_s)\d s+2\E\int_0^{t\land \tau_1} |\nn f|^2(X_s)\d
s.\end{split}\end{equation*} So,

\beg{equation}\label{R1}\beg{split}  &\E f^2(X_{t\land \tau_1})-\{\E
f(X_{t\land \tau_1})\}^2  = 2\int_0^{t\land \tau_1}
\{f(X_s)-f(X_0)\}Lf(X_s)\d s\\
&\qquad -\bigg(\E\int_0^{t\land \tau_1} Lf(X_s)\d s \bigg)^2+
2\E\int_0^{t\land \tau_1} |\nn f|^2(X_s)\d
s.\end{split}\end{equation} Since $Lf$ is bounded on $\B_1:=\{x:
\rr(o,x)\le 1\},$ we have

\beq\label{R2}\bigg(\E\int_0^{t\land \tau_1} Lf(X_s)\d s \bigg)^2\le
c t^2\end{equation}  for some $c>0.$ Moreover, due to Proposition
\ref{P1},

\beq\label{R0} \P(\tau_1\le t)\le c_1\e^{-c_2/t},\ \ \
t>0\end{equation} holds for some constants $c_1,c_2>0.$ Thus,

\beg{equation}\label{R3}\beg{split}&\big| P_tf^2(o) -(P_tf)^2(o)
-\big(\E f^2(X_{t\land \tau_1})-\{\E f(X_{t\land
\tau_1})\}^2\big)\big|=\texto(t^2),\\
& \E\int_0^{t\land \tau_1} |\nn f|^2(X_s)\d s = t|\nn
f(o)|^2+\int_0^t \E \big\{|\nn f|^2(X_{s\land \tau_1})-|\nn
f(o)|^2\big\}\d s+\texto(t^2),\end{split}\end{equation} where and in
what follows, $\texto(s)$ stands for a function of $s>0$ such that
$\lim_{s\to 0} \texto(s)/s =0.$

Similarly,  applying the It\^o formula to $\{f(X_s)-f(o)\}
Lf(X_{s})$, we obtain (note that $Nf|_{\pp M}=0$)

\beg{equation}\label{**1}\beg{split} &\E\int_0^{t\land \tau_1}
\{f(X_s)-f(o)\}Lf(X_s)\d s\\
&= \texto(t^2) +\int_0^t
\E\big[(f(X_{s\land\tau_1})-f(o))Lf(X_{s\land \tau_1})
\big]\d s\\
&= \texto(t^2) + \E\int_0^{t}\d s\int_0^{s\land\tau_1}
L\{(f-f(o))Lf\} (X_r)\d r \\
&\qquad +   \E\int_0^t \d s\int_0^{s\land\tau_1} \{(f-f(o))N
Lf\}(X_r)\d l_r.\end{split}\end{equation} Noting that

$$f(X_r)-f(o)= \ss 2 \int_0^r \<\nn f(X_u), \Phi_u\circ \d B_u\>
+\int_0^r Lf(X_u)\d u,\ \ \ u\le \tau_1,$$  and that

$$\E \sup_{r\in [0,t]}\bigg(\int_0^r \<\nn f(X_u), \Phi_u\circ\d B_u\>\bigg)^2
\le c_2t,\ \ t\in [0,1]$$ holds for some constant $c_2>0$, we obtain
from (\ref{**1}) and (\ref{A4}) that

\beq\label{R4} \bigg|\E\int_0^{t\land \tau_1}
\{f(X_s)-f(o)\}Lf(X_s)\d s\bigg|\le c_3t^2,\ \ t\in
[0,1]\end{equation} holds  for some constant $c_3>0.$ Finally, by
Theorem \ref{T}(1), we have

\beq\label{R5}\E|\nn f|^2(X_{s\land \tau_1})= |\nn f|^2(o) +
\ff{4\ss t}{\ss\pi} \II(\nn f,\nn f)(o)+
\texto(t^{1/2})\end{equation} for small $t>0.$ Combining this with
(\ref{R1}), (\ref{R2}), (\ref{R3})  and (\ref{R4}), and noting that
$U=\nn f(o)$, we conclude that

\beq\label{R6}P_tf^2(o) -(P_tf)^2(o)=2t |\nn f(o)|^2 +
\ff{16t^{3/2}}{3\ss\pi}\II(X,X)+\texto(t^{3/2}).\end{equation}
Finally, (\ref{R5}) and (\ref{R0}) imply that

$$\ff{\e^{2Kt}-1}K P_t|\nn f|^2(o) = 2t|\nn f(o)|^2 +\ff{8t^{3/2}}{\ss\pi}\II(U,U)+\texto(t^{3/2}). $$
Since $\ff {16} 3<8,$ combining this with (10) and (\ref{R6}) we
conclude that $\II(U,U)\ge 0.$

(b) Let $X_0=o\in M\setminus \pp M$, we aim to show that $ \Ric-\nn
Z\ge -K$ holds on $T_oM.$  Let $R>0$ such that $\B_R\cap \pp
M=\emptyset.$ Since $l_t$ increases only when $X_t\in \pp M$,
$l_t=0$ for $t\le \tau_R.$ Hence, due to Proposition \ref{P1}, for
any $f\in C_b^\infty(M)$,

 \beg{equation}\label{AB1}\beg{split} &P_t f^2(o)-(P_tf)^2(o) = \texto(t^2) + \E f^2(X_{t\land \tau_R})
 -\big(\E f(X_{t\land \tau_R})\big)^2 \\
 &  = \texto(t^2) + \int_0^t\big\{\E Lf^2(X_{s\land\tau_R})  -2f(o) \E Lf(X_{s\land\tau_R}) \big\}\d s
 -\bigg(\int_0^t\E Lf(X_{s\land\tau_R})\d
 s\bigg)^2.\end{split}\end{equation}
By the continuity of $s\mapsto Lf(X_{s\land\tau_R})$, we have

\beq\label{AB2} \bigg(\int_0^t\E Lf(X_{s\land\tau_R})\d s\bigg)^2=
(Lf)^2(o) t^2 +\texto(t^2).\end{equation} Similarly, it is easy to
see that

\beg{equation*}\beg{split} &\E  Lf^2(X_{s\land\tau_R})  -2f(o) \E
Lf(X_{s\land\tau_R})\\
&=Lf^2(o)-2f(o)Lf(o)+s\big\{LLf^2-2fLLf\big\}(o)+\texto(s)\\
& = 2|\nn f|^2(o) + 2s \{L|\nn f|^2(o) +(Lf)^2(o)+2\<\nn f,\nn
Lf\>(o)\}+ \texto(s).\end{split}\end{equation*} Combining this with
(\ref{AB1}) and (\ref{AB2}) we obtain

\beq\label{AB3}P_t f^2(o)-(P_tf)^2(o)= 2t|\nn f|^2(o) +t^2 (L|\nn
f|^2+2\<\nn f,\nn Lf\>\}(o)+ \texto(t^2).\end{equation} Finally, by
Proposition \ref{P1} and noting that $l_s=0$ for $s\le \tau_R$, we
have

$$P_t|\nn f|^2(o) = \texto(t^2) + \E |\nn f|^2(X_{t\land\tau_R}) =
|\nn f|^2(o) + t L |\nn f|^2(o) + \texto(t).$$ Combining this with
(10) and (\ref{AB3}), we conclude that

$$\ff 1 2 L|\nn f|^2(o)-\<\nn f, \nn Lf\>(o)\ge -K|\nn f|(o),\ \ f\in C_b^\infty(M).$$ This
completes the proof by the Bochner-Weitzenb\"ock formula.\end{proof}

\section{Proof of Theorem \ref{T1.1}}

 By taking $\mu=\dd_o$, we have $\mu_F^T= \Pi_o^T(F)\dd_o=\dd_o$. So,
(3) follows from each of (2), (7) and (8).    Next, (4) follows from
(3) by taking $F(X_{[0,T]}) = f(X_T),$ and (5) implies  (6) by
taking $p=2$ and $\mu=\dd_x,\nu=\dd_y$. Moreover, it is clear that
(8) follows from (7) while (7) is implied by (2) and (5). So, it
suffices to prove that $(1) \Rightarrow (3)\Rightarrow (2), (4)
\Rightarrow (1) \Rightarrow (6)\Rightarrow (5)$ and $(6)\Rightarrow
(1)$, where $``\Rightarrow$" stands for $``$implies".

(a) $(1) \Rightarrow (3)$.  We shall only consider the case where
$\pp M$ is non-empty and convex. For the case without boundary, the
following argument works well by taking $l_t=0$ and $N=0$. The idea
of the proof comes from \cite{WZ}, where elliptic diffusions on
$\R^d$ were concerned. Let $B_t$ be the $d$-dimensional Brownian
motion on the naturally filtered probability space $(\OO, \F_t,
\P).$
  Let $\{X_t:\ t\ge 0\}$ solve (\ref{2.0}) with $X_0=o$.

Next, let $F$ be a positive bounded measurable function on $M^T$
such that $\inf F>0$ and $\Pi_o^T(F)=1.$ Then

$$m_t:= \E_\P(F(X_{[0,T]})|  \F_t)\ \ \text{and}\ \ L_t :=\int_0^t \ff{\d m_s}{m_s},\ \ \ t\in [0,T]$$ are
square-integrable $\F_t$-martingales under $\P,$ where $\E_\P$ is
the expectation taken for the probability measure $\P.$ Obviously,
we have

\beq\label{2.1} m_t= \e^{L_t-\ff 1 2 \<L\>_t},\ \ \ t\in
[0,T].\end{equation} Since $\scr F_t$ is the natural filtration of
$B_t$, by the martingale representation theorem (cf. \cite[Theorem
6.6]{IW}), there exists a unique $\F_t$-predictable process $\bb_t$
on $\R^d$ such that

\beq\label{2.2} L_t= \int_0^t \<\bb_s,\d B_s\>,\ \ \ t\in
[0,T].\end{equation} Let $\d\Q= F(X_{[0,T]})\d \P.$ Since $\E_\P
F(X_{[0,T]})=\Pi_\mu^T(F)=1,$
 $\Q$ is a probability measure on $\OO$. Due to (\ref{2.1}) and
 (\ref{2.2}) we have

 $$F(X_{[0,T]})=m_T= \e^{\int_0^T \<\bb_s,\d B_s\>-\ff 1 2 \int_0^T\|\bb_s\|^2\d s}.$$
 Moreover, by
 the Girsanov theorem,

\beq\label{2.3}\tt B_t:= B_t-\int_0^t \bb_s\d s,\ \ \ t\in
[0,T]\end{equation} is a $d$-dimensional Brownian motion under the
probability measure $\Q$.

Let $Y_t$ solve the SDE

\beq\label{2.4} \d Y_t= \ss 2\,P_{X_t, Y_t} \Phi_t\circ \d \tt B_t+
Z(Y_t)\d t-N(Y_t)\d\tt l_t,\ \ \ Y_0=o,\end{equation} where $P_{X_t,
Y_t}$ is the parallel displacement along the minimal geodesic from
$X_t$ to $Y_t$ and $\tt l_t$ is the local time of $Y_t$ on $\pp M$.
As explained in e.g. \cite[Section 3]{ATW06}, we may assume that the
minimal geodesic is unique so that $P_{x,y}$ is smooth in $x,y\in
M$. Since, under $\mathbb Q$,  $\tt B_t$ is a $d$-dimensional
Brownian motion, the distribution of $Y_{[0,T]}$ is $\Pi_{o}^T$.

On the other hand, by (\ref{2.0}) and (\ref{2.3}), we have

\beq\label{2.5} \d X_t= \ss 2\,\Phi_t\circ \d \tt B_t +  Z(X_t) +\ss
2\,\Phi_t \bb_t \d t- N(X_t)\d l_t.\end{equation} Since for any
bounded measurable function $G$ on $M^T$

$$\E_\Q G(X_{[0,T]})= \E_\P(FG)(X_{[0,T]}) =\Pi_\mu^T(FG),$$ we conclude that under
$\Q$ the distribution of $X_{[0,T]}$ coincides with $F\Pi_\mu^T$.
Therefore,

\beq\label{2.6} W_{2,\rr_\infty}(F\Pi_\mu^T,\Pi_{\mu_F^T}^T)^2 \le
\E_\Q \rr_\infty(X_{[0,T]}, Y_{[0,T]})^2 =\E_\Q\max_{t\in [0,T]}
\rr(X_t, Y_t)^2.\end{equation} By   the convexity of $\pp M$ we have

$$\<N(x),\nn\rr(y,\cdot)(x)\>= \<N(x), \nn \rr(\cdot,y)(x)\>\le 0,\ \ \ x\in \pp M.$$
 Combining this with the It\^o formula for $(X_t,Y_t)$ given by (\ref{2.4}) and (\ref{2.5}),
we obtain from (\ref{C}) that

\beg{equation*}\beg{split} \d\rr(X_t,Y_t) & \le   K\rr(X_t,Y_t) \d t
+\ss 2\,
\<\Phi_t\bb_t, \nn \rr(\cdot, Y_t)(X_t)\>\d t\\
&\le \Big( K\rr(X_t,Y_t) +\ss 2\,\|\bb_t\|\Big)\d
t,\end{split}\end{equation*}  see e.g. \cite[Lemmas 2.1 and
2.2]{W94}. Since we are using the coupling by parallel displacement
instead of the mirror reflection, the martingale part here
disappears (cf. Theorem 2 and (2.5) in \cite{K}). Since $X_0=Y_0,$
this implies

$$\rr(X_t,Y_t)^2 \le \e^{2Kt} \bigg(\ss 2\,\int_0^t \e^{-Ks}\|\bb_s\|\, \d s\bigg)^2\le
\ff{\e^{2Kt}-1}K\int_0^t \|\bb_s\|^2\d s,\ \ \ t\in [0,T].$$
Therefore,

\beq\label{2.7}  \E_\Q\max_{t\in [0,T]}\rr(X_t,Y_t)^2 \le
\ff{\e^{2KT}-1}K \int_0^T \E_\Q\|\bb_s\|^2\d s.\end{equation} It is
clear that

\beq\label{2.7'} \beg{split} &\mathbb E_{\mathbb Q} \|\bb_s\|^2 =
\E_\P\big(m_T\|\bb_s\|^2\big)\\
&= \E_\P\big( \|\bb_s\|^2\E_\P(m_T|\F_s)\big) =
\E_\P\big(m_s\|\bb_s\|^2\big),\ \ s\in
[0,T].\end{split}\end{equation} Finally, since (\ref{2.1}) and
(\ref{2.2})  yield

$$\d\<m\>_t = m_t^2 \d\<L\>_t=m_t^2 \|\bb_t\|^2 \d t,$$ we have

\beg{equation*}\beg{split} \d m_t\log m_t &= (1+\log m_t)\d m_t +\ff{\d\<m\>_t} {2m_t}\\
&= (1+\log m_t) \d m_t +\ff {m_t} 2 \|\bb_t\|^2 \d
t.\end{split}\end{equation*} As $m_t$ is a $\P$-martingale,
combining this with (\ref{2.7'}) we obtain

\beq\label{2.7''} \int_0^T \mathbb E_{\mathbb Q} \|\bb_s\|^2\d s =2
\E_\P F(X_{[0,T]})\log F(X_{[0,T]}).\end{equation} Therefore,
(\ref{1.3}) follows from (\ref{2.6}), (\ref{2.7}) and (\ref{2.7''}).

(b) $(3)\Rightarrow (2)$. By (3), for each $x\in M$, there exists

$$\pi_x\in \scr C\Big(\ff{F}{\Pi_x^T(F)}\Pi_x^T, \Pi_x^T\Big)$$ such
that

\beq\label{NM}\int_{M^T\times M^T}
\rr_\infty(\gg,\eta)^2\pi_x(\d\gg,\d\eta) \le \ff 2 K(\e^{2KT}-1)
\Pi_x^T\Big(\ff{F}{\Pi_x^T(F)}\log \ff F
{\Pi_x^T(F)}\Big).\end{equation} If $x\mapsto \pi_x(G)$ is
measurable for bounded continuous functions $G$ on $M^T\times M^T$,
then

$$\pi:=\int_M \pi_x\mu_F^T(\d x)\in \scr C(F\Pi_\mu^T,
\Pi_{\mu_F^T}^T)$$ is well defined and by (\ref{NM})

\beg{equation*}\beg{split} \int_{M^T\times M^T} \rr_\infty^2\d\pi
&\le \ff 2 K(\e^{2KT}-1) \Pi_x^T\Big(F\log \ff F
{\Pi_x^T(F)}\Big)\mu(\d x) \\
&\le \ff 2 K(\e^{2KT}-1)\Pi_\mu^T(F\log
F).\end{split}\end{equation*} This implies the inequality in (2).

To confirm the measurability of $x\mapsto \pi_x$, we first consider
discrete $\mu$, i.e. $\mu= \sum_{n=1}^\infty \vv_n \dd_{x_n}$ for
some $\{x_n\}\subset M$ and $\vv_n\ge 0$ with $\sum_{n=1}^\infty
\vv_n=1.$ In this case

$$\pi_x= \sum_{n=1}^\infty 1_{\{x=x_n\}}\pi_{x_n},\ \
\mu\text{-a.e.}$$  which is measurable in $x$ and $\pi=
\sum_{n=1}^\infty \mu_F^T(\{x_n\}) \pi_{x_n}.$ Hence, the inequality
in (2) holds. Then,  for general $\mu$, the desired inequality can
be derived by approximating $\mu$ with discrete distributions in a
standard way, see (b) in the proof of \cite[Theorem 4.1]{FWW}.

(c) $(4)\Rightarrow (1)$. According to \cite[Section 7]{OV} (see
also \cite[Section 4.1]{BGL}), by first applying the
transportation-cost inequality in (3) to $1-\vv+\vv f$ in place of
$f$, then letting $\vv\to 0$, we obtain the Poincar\'e inequality

\beq\label{P} P_Tf^2\le \ff {\e^{2KT}-1}K P_T|\nn f|^2 +(P_Tf)^2,\ \
\ f\in C_b^1(M), T>0.\end{equation}  Thus, the proof is finished by
Theorem \ref{T2.4}.

(d) $(1) \Rightarrow (6)$.  Let $X_t$ solve (\ref{2.0}) with $X_0=x$
and  $Y_t$ solve

\beq\label{4.1} \d Y_t= \ss 2\,P_{X_t, Y_t} \Phi_t\circ \d B_t+
Z(Y_t)\d t-N(Y_t)\d\tt l_t, \ Y_0=y,\end{equation} where $\tt l_t$
is the local time of $Y_t$ on $\pp M$.  Since $\pp M$ is convex and
$(\ref{C})$ holds, as explained in (a), we have

$$\d\rr(X_t,Y_t)\le K \rr(X_t,Y_t)\d t.$$
Thus, $\rr_\infty(X_\cdot, Y_\cdot)\le \e^{KT}\rr(x,y)$. This
implies (6).

(e) $(6) \Rightarrow (5)$. By (6), for any $x,y\in M$, there exists
$\pi_{x,y}\in \scr C(\Pi_x^T,\Pi_y^T)$ such that

$$\int_{M^T\times M^T}\rr_\infty^p \d \pi_{x,y}\le
\e^{KT}\rr(x,y)^p.$$ As explained in (b),  we  assume that $\mu$ and
$\nu$ are discrete, so that for any $\pi^0\in \scr (\mu,\nu)$,
$\pi_{x,y}$ has a $\pi^0$-version measurable in $(x,y)$. Thus,

$$\pi:= \int_{M\times M} \pi_{x,y}\pi^0(\d x,\d y)\in \scr
C(\Pi_\mu^T, \Pi_\nu^T)$$ satisfies

$$\int_{M^T\times M^T} \rr_\infty^p\d\pi \le \e^{KT} \int_{M\times
M}\rr(x,y)^p \pi^0(\d x,\d y).$$ This implies the desired inequality
in (5).

(f) $(6) \Rightarrow (1)$. Let $T>0$ be fixed. For any $x,y\in M$,
let $\pi_{x,y}\in \scr C(P_T(x,\cdot), P_T(y,\cdot))$ be the optimal
coupling for $W_{2,\rr}$, i.e.

\beq\label{LLL}W_{2,\rr}(P_T(x,\cdot), P_T(y,\cdot))^2
=\int_{M\times M}\rr^2\d\pi_{x,y}.\end{equation} Then for any $f\in
C_b^2(M)$,  (6) implies

\beg{equation}\label{4.2} \beg{split}
&\ff{|P_Tf(x)-P_Tf(y)|}{\rr(x,y)}\le \int_{M\times M}
\ff{|f(z_1)-f(z_2)|}{\rr(z_1,z_2)}\cdot
\ff{\rr(z_1,z_2)}{\rr(x,y)}\pi_{x,y}(\d z_1,\d z_2)\\
&\le\ff{ W_{2,\rr}(P_T(x,\cdot),
P_T(y,\cdot))}{\rr(x,y)}\bigg\{\int_{M\times M}
\ff{(f(z_1)-f(z_2))^2}{\rr(z_1,z_2)^2}\pi_{x,y}(\d z_1,\d
z_2)\bigg\}^{1/2}\\
&\le \e^{KT} \bigg\{\int_{M\times M}
\ff{(f(z_1)-f(z_2))^2}{\rr(z_1,z_2)^2}\pi_{x,y}(\d z_1,\d
z_2)\bigg\}^{1/2}.\end{split}\end{equation}Noting that $f\in
C_b^2(M)$ implies

$$|f(z_1)-f(z_2)|^2\le \rr(z_1,z_2)^2 |\nn f|^2(z_1) + c
\rr(z_1,z_2)^3$$ for some constant $c>0$, by (6) and (\ref{LLL}) we
obtain

$$\int_{M\times M}
\ff{(f(z_1)-f(z_2))^2}{\rr(z_1,z_2)^2}\pi_{x,y}(\d z_1,\d z_2)\le
P_T|\nn f|^2(x) +c\e^{KT}\rr(x,y).$$ Therefore,  letting $y\to x$ in
(\ref{4.2}) we arrive at

$$|\nn P_T f(x)| \le \e^{KT} (P_T|\nn f|^2(x))^{1/2}.$$ By a standard
argument of Bakry and Emery, this implies the Poincar\'e inequality
(\ref{P}). Thus, (1) holds according to Theorem \ref{T2.4}.

\section{The case with  diffusion coefficient}

Let $\psi>0$ be a smooth function on $M$, and let $\Pi_{\mu,\psi}^T$
be the distribution of the $($reflecting if $\pp M\ne\emptyset)$
diffusion process generated by $L_\psi:= \psi^2(\DD+Z)$ on time
interval $[0,T]$ with initial distribution $\mu$, and let
$\Pi_{x,\psi}^T=\Pi_{\dd_x,\psi}^T$ for $x\in M$. Moreover, for
$F\ge 0$ with $\Pi_{\mu,\psi}^T(F)=0$, let

$$\mu_{F,\psi}^T(\d x)= \Pi_{x,\psi}^T(F)\mu(\d x).$$

\beg{thm}\label{T5.1} Assume that $\pp M$ is either empty or convex
and let  $(\ref{C})$ hold. Let $\psi\in C_b^\infty(M)$ be strictly
positive. Let

$$K_\psi= K^+\|\psi\|_\infty^2 + 2\|Z\|_\infty \|\nn\psi\|_\infty\|\psi\|_\infty.$$
Then

$$W_{2,\rr_\infty}(F\Pi_{\mu,\psi}^T, \Pi_{\mu_{F,\psi}^T,\psi}^T)^2\le 2 C(T,\psi)\Pi_{\mu,\psi}^T(F\log F),\
\ \ \mu\in\scr P(M), F\ge 0, \Pi_{\mu,\psi}^T(F)=1$$ holds for

$$C(T,\psi):= \inf_{R>0} \Big\{(1+R^{-1}) \|\psi\|_\infty^2 \ff{\e^{2K_\psi T}-1}{K_\psi}
\exp\Big[2(1+R)\|\nn \psi\|_\infty^2 \ff{\e^{2K_\psi
T}-1}{K_\psi}\Big]\Big\}.$$\end{thm}

 \beg{proof}  As explained in (a) of the
proof of Theorem \ref{T1.1}, we shall only consider the case that
$\pp M$ is non-empty and convex. According to the proof of
$``(3)\Rightarrow (2)$", it suffices to prove for $\mu=\dd_o, o\in
M$. In this case the desired inequality reduces to

\beq\label{NM2} W_{2,\rr_\infty}(F\Pi_{o,\psi}^T,\Pi_{o,\psi}^T)\le
2C(T,\psi) \Pi_{o,\psi}^T(F\log F),\ \ F\ge 0,
\Pi_{o,\psi}^T(F)=1.\end{equation}

Since   the diffusion coefficient is non-constant, it is convenient
to adopt the It\^o differential $\d_I$ for the Girsanov
transformation. So, the reflecting diffusion process generated by
$L_\psi:= \psi^2(\DD+Z)$ can be constructed by solving the It\^o SDE

\beq\label{5.0} \d_I X_t= \ss 2\, \psi(X_t) \Phi_t \d B_t
+\psi^2(X_t) Z(X_t)\d t +N(X_t)\d l_t,\end{equation} where   $X_0=o$
and $B_t$ is the $d$-dimensional Brownian motion with natural
filtration $\scr F_t$. Let $\bb_t, \mathbb Q$ and $\tt B_t$ be fixed
in the proof of Theorem \ref{T1.1}. Then

\beq\label{5.1} \d_I X_t= \ss 2\, \psi(X_t) \Phi_t \d \tt B_t
+\big\{\psi^2(X_t) Z(X_t)+\ss 2\, \psi(X_t)\Phi_t\bb_t\big\}\d t
+N(X_t)\d l_t.\end{equation} Let $Y_t$ solve

\beq\label{5.2} \d_I Y_t= \ss 2\, \psi(Y_t) P_{X_t, Y_t}\Phi_t \d
\tt B_t +\psi^2(Y_t) Z(Y_t)\d t +N(Y_t)\d \tt l_t,\ \ \
Y_0=o,\end{equation} where $\tt l_t$ is the local time of $Y_t$ on
$\pp M$. As  in (a) of the proof of Theorem \ref{T1.1}, under
$\mathbb Q$, the distributions of $Y_{[0,T]}$ and $X_{[0,T]}$ are
$\Pi_{o,\psi}^T$ and $F\Pi_{o,\psi}^T$ respectively. So,

\beq\label{T55} W_{2,\rr_\infty}(F\Pi_{o,\psi}^T,
\Pi_{o,\psi}^T)^2\le \mathbb E_{\mathbb Q} \max_{t\in [0,T]}\rr(X_t,
Y_t)^2.\end{equation}  Noting that due to the convexity of $\pp M$

$$\<N(x),\nn \rr(y,\cdot)(x)\>= \<N(x),\nn\rr(\cdot,y)(x)\>\le 0,\ \ x\in\pp M,$$
 by (\ref{5.1}), (\ref{5.2}) and the It\^o
formula,  we obtain

\beg{equation}\label{5.3}\beg{split} \d \rr(X_t,Y_t) \le & \ss
2\,\big\{\psi(X_t) \<\nn \rr(\cdot, Y_t)(X_t), \Phi_t\d\tt B_t\>\\
&\qquad\qquad\qquad+\psi(Y_t) \<\nn
\rr(X_t,\cdot)(Y_t), P_{X_t, Y_t}\Phi_t\d\tt B_t\>\big\}\\
&+\Big\{\sum_{i=1}^{d-1} U_i^2 \rr (X_t, Y_t) +
\<\psi(X_t)^2Z(X_t)+\ss 2\,\psi(X_t) \Phi_t\bb_t , \nn\rr(\cdot,
Y_t)(X_t)\>\\
&\qquad\qquad\qquad +\psi(Y_t)^2
\<Z(Y_t),\nn\rr(X_t,\cdot)(Y_t)\>\Big\}\d
t,\end{split}\end{equation}where $\{U_i\}_{i=1}^{d-1}$ are vector
fields on $M\times M$ such that $\nn U_i (X_t, Y_t)=0$ and

$$U_i(X_t,Y_t) = \psi(X_t)V_i +\psi(Y_t)P_{X_t, Y_t} V_i,\
\ \ 1\le i\le d-1$$ for $\{V_i\}_{i=1}^d$ an OBN of $T_{X_t}M$ with
$V_d= \nn \rr(\cdot, Y_t)(X_t).$

In order to calculate $U_i^2\rr(X_t,Y_t)$, we adopt the second
variational formula for the distance. Let $\rr_t=\rr(X_t,Y_t)$ and
let $\{J_i\}_{i=1}^{d-1}$ be Jacobi fields along the minimal
geodesic $\gg: [0,\rr_t]\to M$ from $X_t$ to $Y_t$ such that
$J_i(0)=\psi(X_t) V_i$ and $J_i(\rr_t)= \psi(Y_t) P_{X_t, Y_t} V_i,
1\le i\le d-1$. Note that the existence of $\gg$ is ensured by the
convexity of $\pp M$. Then, by the second variational formula and
noting that $\nn U_i(X_t,Y_t)=0$, we have

\beg{equation}\label{5.3'} I := \sum_{i=1}^{d-1} U_i^2 \rr(X_t, Y_t)
= \sum_{i=1}^{d-1} \int_0^{\rr_t} \big\{ |\nn_{\dot\gg} J_i|^2
-\<\scr R(\dot \gg, J_i) J_i, \dot \gg\>\big\}(s)\d s,
\end{equation} where $\scr R$ is the curvature tensor. Let

$$\tt J_i(s)= \Big(\ff s {\rr_t} \psi(Y_t) +\ff{\rr_t-s}{\rr_t}
\psi(X_t)\Big)P_{\gg(0), \gg(s)} V_i,\ \ 1\le i\le d-1.$$ We have
$\tt J_i(0)= J_i(0)$ and $\tt J_i(\rr_t)= J_i(\rr_t), 1\le i\le
i-1.$ By the index lemma,

\beg{equation}\label{5.4} \beg{split} I &\le\sum_{i=1}^{d-1}
\int_0^{\rr_t} \big\{ |\nn_{\dot\gg} \tt J_i|^2 -\<\scr R(\dot \gg,
\tt
J_i) \tt J_i, \dot \gg\>\big\}(s)\d s\\
&= \ff{(d-1)(\psi(X_t)-\psi(Y_t))^2}{\rr_t}\\
&\qquad-\ff 1 {\rr_t^2} \int_0^{\rr_t} \big\{s\psi(Y_t)+
(\rr_t-s)\psi(X_t)\big\}^2\Ric\big(\dot\gg(s),\dot\gg(s)\big)\d
s.\end{split}\end{equation} Moreover,

\beg{equation}\label{5.5}\beg{split} &\psi(X_t)^2\<Z(X_t),
\nn\rr(\cdot, Y_t)(X_t)\> +\psi(Y_t)^2
\<Z(Y_t),\nn\rr(X_t,\cdot)(Y_t)\>\\
&=\ff 1 {\rr_t^2} \int_0^{\rr_t} \ff{\d}{\d s}
\Big\{\big(s\psi(Y_t)+(\rr_t-s)\psi(X_t)\big)^2\<Z(\gg(s)),
\dot\gg(s)\>\Big\}\d s\\
&=\ff 1 {\rr_t^2} \int_0^{\rr_t}
 \big(s\psi(Y_t)+(\rr_t-s)\psi(X_t)\big)^2\<(\nn_{\dot\gg} Z)\circ\gg,
\dot\gg\>(s) \d s\\
&\qquad+ \ff 2 {\rr_t^2} \int_0^{\rr_t} \<Z\circ\gg,\dot \gg \>(s)
(\psi(Y_t)-\psi(X_t))\big(s\psi(Y_t)+(\rr_t-s)\psi(X_t)\big)\d
s.\\
&\le \ff 1 {\rr_t^2} \int_0^{\rr_t}
 \big(s\psi(Y_t)+(\rr_t-s)\psi(X_t)\big)^2\<(\nn_{\dot\gg} Z)\circ\gg,
\dot\gg\>(s) \d s+2\|Z\|_\infty \|\psi\|_\infty\|\nn\psi\|_\infty
\rr_t.\end{split}\end{equation} Finally, we have

$$\<\nn\rr (X_t,\cdot)(Y_t), P_{X_t, Y_t}\Phi_t\d\tt
B_t\>= \<P_{Y_t, X_t}\nn \rr(X_t,\cdot)(Y_t),\Phi_t\d\tt
B_t\>=-\<\nn \rr(\cdot,Y_t)(X_t),\Phi_t\d\tt B_t\>.$$ Combining this
with (\ref{5.3}), (\ref{5.3'}), (\ref{5.4}) and (\ref{5.5}),  we
arrive at

\beg{equation}\label{5.6}\beg{split} \d\rr(X_t,Y_t) \le &\ss 2\,
(\psi(X_t)-\psi(Y_t))\<\nn\rr(\cdot, Y_t)(X_t),\Phi_t\d\tt B_t\>\\
&+ K_\psi \rr(X_t,Y_t)\d t  + \ss 2\, \|\psi\|_\infty \|\bb_t\|\d
t=:\d N_t.\end{split}\end{equation} Then

$$M_t:= \ss 2 \, \int_0^t\e^{-K_\psi s}(\psi(X_s)-\psi(Y_s))\<\nn\rr(\cdot, Y_s)(X_s),\Phi_s\d\tt
B_s\>$$ is a $\mathbb Q$-martingale such that

\beq\label{5.6'}\rr(X_t,Y_t)\le   \e^{K_\psi t}M_t + \ss 2 \,
\e^{K_\psi t} \int_0^t \e^{-K_\psi s}\|\psi\|_\infty\|\bb_s\|\d s ,\
\ \ t\in [0,T].\end{equation}  So, by the Doob inequality we obtain

\beg{equation*}\beg{split} h_t &:= \mathbb E_{\mathbb Q}\max_{s\in
[0,t]}
\rr(X_s,Y_s)^2\\
&\le (1+R)\e^{2K_\psi t} \mathbb E_{\mathbb Q} \max_{s\in [0,t]}
M_s^2\d s+ 2\|\psi\|_\infty^2(1+R^{-1}) \e^{2K_\psi
t}\mathbb E_{\mathbb Q}\bigg(\int_0^t \e^{-K_\psi s}\|\bb _s\|\d s\bigg)^2\\
&\le 4 (1+R) \e^{2K_\psi t}\,\mathbb E_{\mathbb Q} M_t^2 +
(1+R^{-1})\|\psi\|_\infty^2 \ff{\e^{2K_\psi t}-1}{K_\psi} \int_0^t
\mathbb E_{\mathbb Q}\|\bb_s\|^2\d s\\
&\le 4 (1+R)\|\nn\psi\|_\infty^2 \e^{2K_\psi t}\int_0^t\e^{-2K_\psi
s}h_s\d s + (1+R^{-1})\|\psi\|_\infty^2 \ff{\e^{2K_\psi
T}-1}{K_\psi} \int_0^t\mathbb E_{\mathbb Q} \|\bb_s\|^2\d
s\end{split}\end{equation*}for any $R>0.$ Since $\e^{-2K_\psi s}$ is
decreasing in $s$ while $h_s$ is increasing in $s$, by the FKG
inequality we have

$$\int_0^t \e^{-2K_\psi s}h_s\d s \le \ff{1-\e^{-2K_\psi t}}{2K_\psi
t} \int_0^t h_s\d s.$$ Therefore,

$$h_t\le 2(1+R)\|\nn\psi\|_\infty^2 \ff{\e^{2K_\psi T}-1}{K_\psi T}
\int_0^t h_s\d s +(1+R^{-1})\|\psi\|_\infty^2 \ff{\e^{2K_\psi
T}-1}{K_\psi} \int_0^t \mathbb E_{\mathbb Q}\|\bb_s\|^2\d s$$ holds
for $t\in [0,T].$ Since $h_0=0,$ this implies that

\beg{equation*} \beg{split} & \mathbb E_{\mathbb Q}\max_{t\in
[0,T]}\rr(X_t,Y_t)^2=h_T\\
&\le (1+R^{-1})\|\psi\|_\infty^2\ff{\e^{2K_\psi T}-1}{K_\psi}
\exp\Big[2(1+R)\|\nn\psi\|_\infty^2\ff{\e^{2K_\psi
T}-1}{K_\psi}\Big] \int_0^T\mathbb E_{\mathbb Q}\|\bb_s\|^2\d
s.\end{split}\end{equation*} Combining this with the (\ref{T55}) and
(\ref{2.7''}),
  we complete the proof.\end{proof}

\beg{thm}\label{T5.2} In the situation of Theorem \ref{T5.1},

$$W_{2,\rr_\infty}(\Pi_{\mu,\psi}^T, \Pi_{\nu,\psi}^T)\le
2\e^{(K_\psi+\|\nn\psi\|_\infty^2)T} W_{2,\rr}(\mu,\nu),\ \ \
\mu,\nu\in \scr P(M), T>0.$$\end{thm}

\beg{proof}  As explained in the proof of $`` (6) \Rightarrow (5)$",
we only consider $\mu=\dd_x$ and $\nu=\dd_y$. Let $X_t$ solve
(\ref{5.0}) with $X_0=x$, and  let $Y_t$ solve, instead of
(\ref{5.2}),

$$\d_I Y_t= \ss 2\, \psi(Y_t) P_{X_t, Y_t}\Phi_t \d
\tt B_t +\psi^2(Y_t) Z(Y_t)\d t +N(Y_t)\d \tt l_t,\ \ Y_0=y.$$ Then,
repeating the proof of Theorem \ref{T5.1}, we have, instead of
(\ref{5.6'}),

\beq\label{5.A2}\rr(X_t,Y_t)\le \e^{K_\psi t} (M_t+ \rr(x,y)),\ \ \
t\ge 0\end{equation}  for

$$M_t:= \ss 2 \, \int_0^t\e^{-K_\psi s}(\psi(X_s)-\psi(Y_s))\<\nn\rr(\cdot, Y_s)(X_s),\Phi_s\d
B_s\>.$$ So,

$$\E\rr(X_t,Y_t)^2\le \e^{2K_\psi t}\bigg\{\rr(x,y)^2
+2\|\nn\psi\|_\infty^2\int_0^t \e^{-2K_\psi s}\E\rr(X_s,Y_s)^2\d
s\bigg\},$$ which implies

$$\E\rr(X_t,Y_t)^2\le \e^{2(K_\psi +\|\nn \psi\|_\infty^2)t}
\rr(x,y)^2.$$ Combining this with   (\ref{5.A2}) and   the Doob
inequality, we arrive at

\beg{equation*}\beg{split} &W_{2,\rr_\infty}(\Pi_{x,\psi}^T,
\Pi_{y,\psi}^T)^2 \le \E \max_{t\in [0,T]} \rr(X_t,Y_t)^2\le
\e^{2K_\psi T} \E\max_{t\in [0,T]} (M_t+\rr(x,y))^2\\
&\le 4    \e^{2K_\psi T} \E(M_T+\rr(x,y))^2= 4\e^{2K_\psi T}
\big(\E M_T^2+\rr(x,y)^2\big)\\
&=4    \e^{2K_\psi T} \bigg(\rr(x,y)^2 + 2\|\nn \psi\|_\infty^2
\int_0^T \e^{-2K_\psi t}\E \rr(X_t,Y_t)^2\d
t\bigg)\\
&\le 4 \e^{2(K_\psi +\|\nn
\psi\|_\infty^2)T}\rr(x,y)^2.\end{split}\end{equation*} This implies
the desired inequality for $\mu=\dd_x$ and $\nu=\dd_y$.
\end{proof}

\section{Extensions to non-convex manifolds}

As explained in the end of Section 1, combining Theorem \ref{T5.1}
with a proper conformal change of metric, we are able to establish
the following transportation-cost inequality on a class of manifolds
with non-convex boundary.

\beg{thm}\label{T6.1} Let $\pp M\ne\emptyset $ with $\II\ge -\si$
for some constant $\si> 0,$ and let $(\ref{C})$ hold for some $K\in
\R.$ Then for any $f\in C_b^\infty(M)$ with $f\ge 1$ and $N\log
f|_{\pp M}\ge \si$, and for any $\mu\in \scr P(M),$

$$ W_{2, \rr_\infty}(F\Pi_{\mu}^T,\Pi_{\mu_F^T}^T)^2\le 2\|f\|_\infty^2 c(T,f)
\Pi_{\mu}^T(F\log F),\ \ F\ge 0, \Pi_{\mu}^T(F)=1$$ holds for

$$c(T,f)=\inf_{R>0} \Big\{(1+R^{-1})   \ff{\e^{2\kk_f T}-1}{\kk_f}
\exp\Big[2(1+R)\|\nn f\|_\infty^2 \ff{\e^{2\kk_f
T}-1}{\kk_f}\Big]\Big\},$$ where

$$ \kk_f= 5\|f\|_\infty\|\nn f\|_\infty\|Z\|_\infty +\big\{2(d-2)
+(d-3)^+\big\} \|\nn f\|_\infty^2 +\|(Kf^2-f\DD f)^+\|_\infty. $$ In
particular,

$$ W_{2, \rr_\infty}(F\Pi_{o}^T,\Pi_o^T)^2\le 2\|f\|_\infty^2 c(T,f)\Pi_o^T(F\log F),\ \ o\in M, F\ge 0,
\Pi_o^T(F)=1.$$

\end{thm}
\beg{proof}   Let $f\in C_b^\infty(M)$ such that $f\ge 1$. Since
$\II\ge -\si$ and $N\log f|_{\pp M}\ge \si$, by \cite[Lemma
2.1]{W07} the boundary $\pp M$ is convex under the new metric

$$\<\cdot,\cdot\>'= f^{-2}\<\cdot,\cdot\>.$$ Let $\DD'$ and $\nn'$ be induced by the new metric.
 Then (see formula (2.2) in \cite{TW})

$$ L = f^{-2}(\DD'+Z'),\ \ \ Z':= f^2 Z + \ff{d-2}2 \nn f^2.$$
 Let $\Ric'$ be the Ricci curvature induced by the new metric, we have (cf. formula
 (3.2) in \cite{FWW})

\beq\label{6.1} \Ric' =\Ric +(d-2)f^{-1}\Hess_f + (f^{-1} \DD f
-(d-3)|\nn \log f|^2)\<\cdot,\cdot\>.   \end{equation} Since the
Levi-Civita connection induced by $\<\cdot,\cdot\>'$ satisfies (cf.
\cite[Theorem 1.59(a)]{B})

$$\nn'_UV = \nn_UV -\<U,\nn\log f\>V -\<V,\nn\log f\>U +\<U,V\> \nn \log f,\ \ \
U,V\in TM,$$ we have

\beg{equation*}\beg{split} &\<\nn'_U Z',U\>' = f^{-2} \big\{\<\nn_U Z', U\> - \<Z',\nn\log f\>|U|^2 \big\}\\
&= 2\<U,\nn\log f\>\<Z,U\> +\<\nn_U Z, U\> + \ff {d-2}{2 f^2}
\Hess_{f^2} (U,U)\\
&\qquad\qquad  -\<Z,\nn\log f\>|U|^2 -\ff{d-2}
2 \<\nn\log f^2, \nn \log f\>|U|^2\\
&\le \<\nn_UZ, U\> + 3 |\nn\log f|\cdot |Z|\cdot |U|^2 +(d-2) f^{-1}
\Hess_f(U,U).\end{split}\end{equation*} Combining this with
(\ref{6.1}), we obtain

    \beg{equation*}\beg{split} &\Ric'(U,U)-\<\nn_UZ', U\>'\\
     &\ge \Ric(U,U)- \<\nn_UZ,U\>
    +\big\{f^{-1}\DD f -(d-3)|\nn\log f| - 3|Z|\cdot |\nn\log f|\big\} |U|^2\\
    &\ge -K'\<U,U\>',\ \ \ U\in TM,\end{split}\end{equation*}
where

\beq\label{KK} K'  = \sup_M\{Kf^2 -f\DD f +(d-3)|\nn f|^2 + 3|Z|f
|\nn f|\}.\end{equation} Noting that   $f\ge 1$, we have

\beg{equation}\label{KK2}\beg{split} & \ss{\<Z',Z'\>'}= f^{-1}
|f^2Z+(d-2)f\nn
f| \le \|f\|_\infty \|Z\|_\infty +(d-2)\|\nn f\|_\infty,\\
&\ss{\<\nn' f^{-1}, \nn' f^{-1}\>'}= f|\nn f^{-1}|\le \|\nn
f\|_\infty.\end{split}\end{equation} Letting $K_\psi$ be defined in
Theorem \ref{T5.1} for the manifold $(M, \<\cdot, \cdot\>')$ and $L=
\psi^2(\DD'+Z')$ with $\psi=f^{-1}$, we deduce from  $f\ge 1$,
(\ref{KK}) and (\ref{KK2}) that

$$K_\psi\le \kk_f.$$ Therefore, $C(T,\psi)\le c(T,f)$ and thus,
Theorem \ref{T5.1} implies

$$W_{2,\rr'_\infty}(F\Pi_{\mu}^T, \Pi_{\mu_F^T}^T)^2\le
2c(T,f)\Pi_\mu^T(F\log F),\ \ \ F\ge 0, \Pi_\mu^T(F)=1,$$  where
$\rr'_\infty$ is the uniform distance on $M^T$ induced by the metric
$\<\cdot,\cdot\>'$. The proof is completed by noting that
$\rr_\infty\le \|f\|_\infty \rr'_\infty.$\end{proof}

Similarly,  since  $K_\psi\le \kappa_f$ and

$$\rr'\le \rr\le \|f\|_\infty \rr,$$  the following result from Theorem \ref{T5.2} by taking $\psi= f^{-1}.$

\beg{thm}\label{T6.2} In the situation of Theorem \ref{T6.1},

$$W_{2,\rr_\infty}(\Pi_\mu^T, \Pi_\nu^T)\le
2\|f\|_\infty\e^{(\kappa_f+\|\nn f^{-1}\|_\infty^2)T}
W_{2,\rr}(\mu,\nu),\ \ \ \mu,\nu\in \scr P(M), T>0.$$\end{thm}

 As a
consequence of Theorems \ref{T6.1} and \ref{T6.2}, we present below
an explicit transportation-cost inequalities for a class of
non-convex manifolds.

\beg{cor}\label{C6.2} Assume that $(\ref{C})$ holds for some $K\ge
0$ and the injectivity radius $\i_{\pp M}$ of $\pp M$ is strictly
positive. Let $\si\ge 0$ and $\gg, k,
>0$ be such that $-\si\le\II\le \gg$ and
$\Sect_M\le k$. Let

$$0<r\le \min\bigg\{\i_{\pp M},\ \ff1 {\ss k} \arcsin\bigg(\ff{\ss
k}{\ss{k+\gg^2}}\bigg)\bigg\}.$$

{\rm (i)} The transportation-cost inequality

$$W_{2, \rr_\infty}(F\Pi_\mu^T,\Pi_{\mu_F^T}^T)^2
\le (2+rd\si)^2 \ff{\e^{2\theta T}-1}{\theta
}\exp\Big[\ff{4(\e^{2\theta T}-1)}{\theta }\Big] \Pi_{\mu}^T(F\log
F)$$ holds for all $\mu\in \scr P(M)$ and $F\ge 0$ with $
\Pi_{\mu}^T(F)=1,$ where

$$ \theta= K\Big(1+rd\si +\ff {r^2d^2\si^2} 4\Big) +\ff {d\si} r \Big(2(d-2)+(d-3)^+ +\ff {d^2}
2\Big)
\si^2 + 5\|Z\|_\infty\si \Big(1+\ff{rd\si}2\Big). $$ In particular,

$$W_{2, \rr_\infty}(F\Pi_o^T,\Pi_o^T)^2
\le (2+rd\si)^2 \ff{\e^{2\theta T}-1}{\theta
}\exp\Big[\ff{4(\e^{2\theta T}-1)}{\theta }\Big] \Pi_o^T(F\log
F)$$holds for all $F\ge 0$ with $ \Pi_\mu^T(F)=1.$

{\rm (ii)} For any $T>0$ and $\mu,\nu\in \scr P(M),$

$$W_{2,\rr_\infty}(\Pi_\mu^T, \Pi_\nu^T)\le
(2+\si r d)\e^{(\theta+\si^2)T} W_{2,\rr}(\mu,\nu).$$
\end{cor}

\beg{proof}   Let

 $$h(s)= \cos\big(\ss k\, s\big)- \ff{\gg } {\ss k}\sin\big(\ss k\, s\big),\
 \ \ s\ge 0.$$ Then $h$ is the unique solution to the equation

 $$h'' + k h=0,\ \ \ h(0)=1, h'(0)= -\gg.$$
 Up to an approximation argument presented in the proof of
\cite[Theorem 1.1]{W05},
  we may apply Theorem \ref{T6.1} to

$$f= 1+ \si \varphi\circ\rr_{\pp M},$$
where $\rr_\pp$ is the Riemannian distance to $\pp M$, which is
smooth on $\{\rr_{\pp M}<\i_{\pp M}\},$ and

 \beg{equation*}\beg{split} &\aa= (1-h(r))^{1-d}\int_0^r
 (h(s)-h(r))^{d-1}\d s,\\
 &\varphi(s) =\ff 1 \aa\int_0^s (h(t)-h(r))^{1-d}\d t\int_{t\land
 r}^r(h(u)-h(r))^{d-1}\d u,\ \ s\ge 0.\end{split}\end{equation*} We
 have $\varphi(0)=1, 0\le \varphi'\le \varphi'(0)=1.$ Moreover, as observed in \cite[Proof of
 Theorem 1.1]{W05},

 $$ \aa\ge \ff r d,\ \  \varphi(r)\le
 \ff{r^2}{2\aa}\le \ff{dr} 2,\ \ \DD\varphi\circ\rr_{\pp M}\ge -\ff 1 {\aa}\ge -\ff d r.$$ So,

\beq\label{LL} \|f\|_\infty\le 1+ \si \varphi(r)\le 1+\ff {r d\si
}2,\ \ \|\nn f\|_\infty\le \varphi'(0)=\si ,\ \ \DD f\ge -\ff {\si
d} r.\end{equation} Noting that (recall that $K\ge 0$)

$$\sup (Kf^2)\le K \Big(1+rd\si  + \ff {r^2 d^2\si^2}4\Big),$$ from (\ref{LL})
we conclude that $\kk_f \le \theta.$ So, (i) follows from
 (\ref{L}) and   \ref{T6.1} for $R=1$, and (ii) follows from
 Theorem \ref{T5.2} and (\ref{LL}).
\end{proof}

\beg{thebibliography}{99}

\bibitem{ATW06}
M. Arnaudon, A. Thalmaier, F.-Y. Wang, \emph{Harnack inequality and
heat kernel
  estimates on manifolds with curvature unbounded below,} Bull. Sci. Math.
  130(2006), 223--233.

\bibitem{ATW09} M. Arnaudon,   A. Thalmaier and F.-Y. Wang,   \emph{Gradient estimates and Harnack inequalities on
non-compact Riemannian manifolds,}  Stoch. Proc. Appl. 2009.

\bibitem{B} A. L. Besse, \emph{Einstein Manifolds,} Springer, Berlin, 1987.

\bibitem{BGL} S. Bobkov, I. Gentil and M. Ledoux, \emph{Hypercontractivity of Hamilton-Jacobi
equations}, J. Math. Pure Appl. 80 (2001), 669--696.

\bibitem{DGW} H. Djellout, A. Guillin, and L.-M. Wu, \emph{Transportation
cost information inequalities and applications to random dynamical
systems and diffusions}, Ann. Probab. 32 (2004),  2702--2732.

\bibitem{FWW} S. Fang, F.-Y. Wang and B. Wu, \emph{Transportation-cost inequality on path spaces with
uniform distance,} Stoch. Proc. Appl. 118(2008), 2181Ð2197.

\bibitem{F} D. Feyel and A. \"Ust\"unel, \emph{Measure
transport on Wiener space and the Girsanov theorem,} C.R. Acad.
Paris 334(2002), 1025--1028.

\bibitem{IW} N. Ikeda and S. Watanabe, \emph{Stochastic Differential Equations,} North-Holland, New York, 1989.


\bibitem{K} W.S. Kendall, \emph{Nonnegative Ricci curvature and the Brownian
coupling property,} Stochastics 19  (1986), 111¨C 129.

\bibitem{L} M. Ledoux, \emph{The geometry of Markov diffusion
generators,} Ann. Facul. Sci. Toulouse 9(2000), 305--366.

\bibitem{RS} M.-K. von Reness and K.-T. Sturm, \emph{Transport
inequalities, gradient estimates, entropy, and Ricci curvature,}
Comm. Pure Math. 58(2005), 923--940.

 \bibitem{OV} F. Otto and C. Villani, \emph{Generalization of an inequality
by Talagrand and links with the logarithmic Sobolev inequality,} J.
Funct. Anal.  173(2000), 361--400.

\bibitem{T} M. Talagrand, \emph{Transportation cost for Gaussian
and other product measures,} Geom. Funct. Anal. 6(1996), 587--600.

\bibitem{TW} A. Thalmaier, F.-Y. Wang, \emph{Gradient estimates for harmonic functions on regular domains in
Riemannian manifolds,} J. Funct. Anal. 155:1(1998),109--124.

 \bibitem{W94} F.-Y. Wang, \emph{Application of coupling methods to the Neumann eigenvalue problem,} Probab.
 Theory Related Fields 98 (1994),  299--306.

 \bibitem{W97} F.-Y. Wang,\emph{On estimation of the logarithmic Sobolev
constant and gradient estimates of heat semigroups,} Probability
Theory Related Fields 108(1997), 87--101.

\bibitem{W02} F.-Y. Wang, \emph{Transportation cost inequalities on path spaces
over Riemannian manifolds,} Illinois J. Math. 46 (2002), 1197--1206.

\bibitem{W04} F.-Y. Wang, \emph{Probability distance inequalities on Riemannian manifolds and path
spaces}, J. Funct. Anal. 206 (2004), 167--190.

\bibitem{W04b} F.-Y. Wang, \emph{Equivalence of dimension-free
Harnack inequality and curvature condition,} Integral Equation and
Operator Theory 48(2004), 547--552.

\bibitem{W05} F.-Y. Wang, \emph{Gradient estimates and the first
Neumann eigenvalue on manifolds with boundary,} Stoch. Proc. Appl.
115(2005), 1475--1486.

\bibitem{W07}  F.-Y. Wang, \emph{Estimates of the first Neumann eigenvalue and the log-Sobolev constant
  on  nonconvex manifolds,} Math. Nachr. 280(2007), 1431--1439.

  \bibitem{W09} F.-Y. Wang, \emph{Second fundamental form and gradient of Neumann semigroups,}  J. Funct. Anal. 256(2009), 3461--3469.

\bibitem{WZ} L.-M. Wu and Z.-L. Zhang, \emph{ Talagrand's $T\sb 2$-transportation inequality
w.r.t. a uniform metric for diffusions}, Acta Math. Appl. Sin. Engl.
Ser. 20 (2004),  357--364.

 \end{thebibliography}

\end{document}